\input amstex
\magnification=\magstep1
\binoppenalty=10000\relpenalty=10000
\documentstyle{amsppt}
\hcorrection{1.3 true cm}
\pagewidth{13.284 true cm}
\pageheight{21.3 true cm}
\topmatter
\title Extraspecial towers and Weil representations\endtitle
\author S. P. Glasby {\rm and} R. B. Howlett\endauthor
\address Dept. of Mathematics, Victoria University of Wellington,
P.O. Box 600, Wellington, New Zealand.\endaddress
\address Department of Pure Mathematics, University of Sydney, NSW 2006,
Australia.\endaddress
\date 30 July 1990\enddate
\endtopmatter

\def\brw{{\bf 1}}
\def\ger{{\bf 2}}
\def\gla{{\bf 3}}
\def\glasby{{\bf 4}}
\def\gor{{\bf 5}}
\def\hup{{\bf 6}}
\def\wall{{\bf 7}}
\def\wardi{{\bf 8}}
\def\wardii{{\bf 9}}

\define\F{\bold F}
\newdimen\timesht\newdimen\timesdp
\def\semiprod{\mathbin{\mathpalette\vbar\times}}
\def\vbar#1#2{\setbox0=\hbox{$#1#2$}\timesht=\ht0\timesdp=\dp0
\advance\timesht -.7pt\advance\timesdp-.2pt
\vrule height\timesht depth\timesdp\mskip-2.8mu\box0}

\define\Aut{\operatorname{Aut}}

\define\gl{\operatorname{GL}}

\define\im{\operatorname{im}}
\hyphenation{ex-tra-spec-ial}

\document

\heading 1. Introduction\endheading

\noindent By an ``extraspecial tower'' we mean a group $G$ with a series of
normal subgroups $G=N_0>N_1>N_2>\cdots$ such that
$N_i\big/N_{i+1}$ is extraspecial and $N_{i+1}$ is contained in the 
derived group of $N_i$ for each $i$. It is the construction of such
groups which is the primary objective of this paper.

The following example sparked our interest in this problem. Let $Q$ and 
$E$ be (respectively) the quaternion group of order~8 and the 
nonabelian group of order 27 and exponent~3. Let $Q^n$ be the central product
of $n$ copies of $Q$ and let $E^n$ be the central product of $n$ copies of $E$,
so that $Q^n$ and $E^n$ are extraspecial groups of orders $2^{2n+1}$ and
$3^{2n+1}$ respectively. The group $\gl_2(\F_3)$, itself a semidirect
product $S_3\semiprod Q$ (where $S_3$ is the symmetric group of degree~3), 
can be embedded in the automorphism group of~$E$. The semidirect product
$(S_3\semiprod Q)\semiprod E$ embeds in $\Aut(Q^3)$, and the semidirect
product $((S_3\semiprod Q)\semiprod E)\semiprod Q^3$ embeds in the
automorphism group of~$E^4$. It is natural to 
ask whether this can be continued indefinitely. In fact it can be, but 
the extraspecial groups start to increase rapidly in size: the next two are
$Q^{81}$ and $E^{2^{80}}$.

The inductive process we use to construct the groups proceeds as 
follows. Let $p_{i-1}$ and $p_i$ be primes, and assume that $G_i$ is a
semidirect product $G_{i-1}\semiprod E_{i-1}$, where $E_{i-1}$ is an
extraspecial $p_{i-1}$-group. An extraspecial $p_i$-group $E_i$ is chosen
so that $E_{i-1}$ acts faithfully on the central quotient
$V_i=E_i\big/Z(E_i)$, regarded as a vector space over a field of
characteristic~$p_i$. We then embed $G_i$ in a group
$A_i$ of automorphisms of $E_i$, so that we may form the semidirect
product $G_{i+1}=G_i\semiprod E_i$ and repeat the process. This embedding is 
accomplished by extending the representation of $E_{i-1}$ on $V_i$ to a
representation of $A_{i-1}\semiprod E_{i-1}$.

The towers which are of most interest to us are ones which are minimal,
in the sense that at each stage the $E_{i-1}$-module $V_i$ is either
irreducible, or as close to irreducible as possible. Without this requirement
there would be no particular difficulty in constructing extraspecial towers.

The groups $A_i$ which arise in the construction are either symplectic 
or unitary groups; the representations of the $A_i\semiprod E_i$ are the 
so-called ``Weil representations'', which have been investigated by many
authors. In our treatment, which is most influenced by Ward [\wardi,\wardii] 
and G\'erardin [\ger], we give explicit splittings of the factor set which
arises in the construction of the Weil representations, including a
determination of the correct sign.\par We would like to thank to referee 
for several helpful comments. 

\heading 2. Basic properties of extraspecial groups\endheading

\noindent In this section we introduce some notation and state, without proof, 
some basic properties of extraspecial groups. Proofs can be found in
Gorenstein [\gor], for instance.

A finite group $P$ of $p$-power order, where $p$ is a prime,
is said to be {\it extraspecial} if $P'$ (the derived group), $Z(P)$
(the centre) and $\Phi(P)$ (the Frattini subgroup) all have order~$p$.
All extraspecial groups are central products of extraspecial groups of 
order~$p^3$, and for each prime~$p$ there are two isomorphism classes of
extraspecial groups of order~$p^3$. For $p=2$ they are represented by the
quaternion group $Q$ and the dihedral group $D$. For $p\ne2$ they are 
represented by a group $E$ of exponent~$p$ and a group $M$ of 
exponent~$p^2$. Every extraspecial group has order $p^{2n+1}$ for some 
integer~$n$, and is isomorphic to precisely one of the central products 
$D^n$ or $D^{n-1}Q$ (if $p=2$), or $E^n$ or $E^{n-1}M$ (if $p\ne2$).
In this paper we will consider extraspecial towers constructed from
the first three of these four types of extraspecial groups.

Let $V$ be a vector space over a field~$K$. In this section
we assume that $K=\F_p$ is a prime field. Let $f\colon V\times V\to K$
be a bilinear form, and let $f^{\roman t}$ be the form defined by
$(x,y)\mapsto f(y,x)$. Assume
that $f-f^{\roman t}$ is nondegenerate; since $f-f^{\roman t}$ is 
alternating, this implies that $\dim_KV$ is even. Furthermore, the set
$E(f)=V\times K$ becomes a group if multiplication is defined by
$$(x_1,z_1)(x_2,z_2)=(x_1+x_2,z_1+z_2+f(x_1,x_2))$$for all $x_1,\,x_2\in 
V$ and $z_1,\,z_2\in K$. Observe that $(0,0)$ is the identity element,
$(x,z)^j=(jx,jz+\tfrac12j(j-1)f(x,x))$ for all integers~$j$, and
$$(x_1,z_1)^{-1}(x_2,z_2)^{-1}(x_1,z_1)(x_2,z_2)=(0,(f-f^{\roman 
t})(x_1,x_2)).$$Since $f-f^{\roman t}$ is nondegenerate we see that both
the centre and derived group are equal to $Z=\{\,(0,z)\mid z\in K\,\}$,
and since $(x,z)^p\in Z$ for all $x$ and~$z$ it follows that $E(f)$
is extraspecial. Furthermore, if $p$ is odd then $(x,z)^p=1$ for all $x$
and~$z$, so that $E(f)$ has exponent~$p$. In summary,
$E(f)\cong E^n$ when $p$ is odd and $\dim V=2n$.

Suppose that $V=\{\,(x,y)\mid x,y\in K\,\}$. In the case $p=2$, let 
$f_D$ be the bilinear form defined by
$$f_D((x_1,y_1),(x_2,y_2))=y_1x_2=\pmatrix x_1&y_1\endpmatrix
\pmatrix0&0\\1&0\endpmatrix\pmatrix x_2\\y_2\endpmatrix$$ and similarly 
let $f_Q$ be the bilinear form with matrix $\left(\smallmatrix 
1&0\\1&1\endsmallmatrix\right)$. It is easily checked that $E(f_D)\cong D$
and $E(f_Q)\cong Q$. In the case when $p$ is an odd prime let $f_E$ be the
bilinear form with matrix
$\left(\smallmatrix 0&1\\-1&0\endsmallmatrix\right)$, so 
that $E(f_E)\cong E$. When $f=f_D$, $f_Q$ or $f_E$ we write elements of 
the group $E(f)$ as triples $(x,y,z)$ of elements of $\F_p$.

Recall that if $f_i\colon V_i\times V_i\to K$ (for $i=1,2$) are bilinear forms
then $f_1\oplus f_2$ defined by
$$(f_1\oplus f_2)((x_1,x_2),(y_1,y_2))=f_1(x_1,y_1)+f_2(x_2,y_2)$$
is a bilinear form on $V_1\oplus V_2$. With $f_D,\,f_Q$ and $f_E$ as 
above it can be seen that the bilinear forms
$\bigoplus_{i=1}^nf_D$, $(\bigoplus_{i=1}^{n-1}f_D)\oplus f_Q$ and 
$\bigoplus_{i=1}^nf_E$ act on $2n$-dimensional spaces and satisfy
$$E(\bigoplus_{i=1}^nf_D)\cong D^n,\quad\hss E((\bigoplus_{i=1}^{n-1}f_D)
\oplus f_Q)\cong D^{n-1}Q,\quad\hss\text{and}\quad
\hss E(\bigoplus_{i=1}^nf_E)\cong E^n.$$

The following well-known proposition can be used to determine the 
isomorphism type of $E(f)$ in the case $p=2$.
\proclaim{2.1 Proposition}Let $f\colon V\times V\to\F_2$ be a bilinear 
form such that $f-f^{\roman t}$ is nondegenerate. The element $(x,z)\in
E(f)$ has order 1 or~2 if and only if $x\in V$ satisfies $f(x,x)=0$. The
number of such $x$ is $2^{n-1}(2^n+1)$ if $E(f)\cong D^n$ and
$2^{n-1}(2^n-1)$ if $E(f)\cong D^{n-1}Q$.\endproclaim
If $g\in\gl(V)$ satisfies $f(x_1g,x_2g)=f(x_1,x_2)$ for all
$x_1,x_2\in V$ then clearly $(x,z)\mapsto(xg,z)$ is an automorphism 
of $E(f)$. Thus, for any group~$G$, a representation $G\to\gl(V)$ which
preserves $f$ gives rise to an action of $G$ on $E(f)$, enabling the 
construction of a semidirect product $G\semiprod E(f)$.
\heading 3. Representations of extraspecial groups\endheading

\noindent In this section we consider faithful irreducible representations
of the groups $E(f)$, and forms they preserve. Parts of this exposition may be
found in \cite{\gor}. Suppose that
$\rho\colon E(f)\to \gl(V')$ is an absolutely irreducible representation, 
where $V'$ is a vector space over a field $K'$. Since the centre of $E(f)$ has
order~$p$ and must be faithfully represented by scalar transformations, $K'$
must contain a primitive $p^{\roman{th}}$ root of~1. Let $\epsilon$ be 
such a root (so that $\epsilon=-1$ if $p=2$).

We are primarily interested in the case when $K'$ is a finite field, of 
characteristic~$p'$, say. Clearly $p'\ne p$. Assuming that $K'$ contains
the necessary primitive $p^{\roman{th}}$ root of~1, we may write down 
faithful absolutely irreducible representations $\rho_E$, $\rho_D$ and 
$\rho_Q$ of $E(f_E)$, $E(f_D)$ and $E(f_Q)$. We define $\rho_E$ by
$$(x,y,z)\mapsto\pmatrix0&1&0&&0\\0&0&1&&0\\&&&\ddots
&\\0&0&0&&1\\1&0&0&&0\endpmatrix^x\pmatrix1&0&&0&0\\0&\epsilon&&0&0\\
&&\ddots&&\\0&0&&\epsilon^{p-2}&0\\0&0&&0&\epsilon^{p-1}\endpmatrix^{2y}
\epsilon^{z-xy}$$ for all $x,y,z\in\F_p$, and $\rho_D$ by
$$(x,y,z)\mapsto\pmatrix0&1\\1&0\endpmatrix^x
\pmatrix1&0\\0&-1\endpmatrix^y(-1)^z$$for all $x,y,z\in\F_2$. Note that 
if $A$ is a matrix satisfying
$A^2=\left(\smallmatrix-1&0\\0&-1\endsmallmatrix\right)$ then 
$(-1)^{x\choose2}A^x$ depends only on the parity of the integer~$x$.
So, abusing notation somewhat, we define $\rho_Q$ by
$$(x,y,z)\mapsto\pmatrix\alpha&\beta\\\beta&-\alpha\endpmatrix^x
\pmatrix0&1\\-1&0\endpmatrix^y(-1)^{z+{x\choose2}+{y\choose2}}
\quad\hbox{for all $x,y,z\in\F_2$}$$
where $\alpha$ and $\beta$ are fixed elements of $K'$ such that
$\alpha^2+\beta^2=-1$. Finiteness of $K'$ guarantees that suitable 
$\alpha$ and $\beta$ always exist.

Let $\rho$ be one of $\rho_D$, $\rho_Q$ or $\rho_E$ as defined above, where 
the $p\times p$ matrices over $K'$ are interpreted as right operators on
the space $V'$ of $p$-component row vectors over~$K'$.
We seek an automorphism $\eta'$ of $K'$ (possibly
the identity) and a $\eta'$-sesquilinear form $f'$ on $V'$
that is preserved by~$\rho$. In matrix terms we have 
$f'(u,v)=uJ(v^{\roman t\eta'})$ for some matrix~$J$, and $\rho$ 
preserves $f'$ if and only if $XJ(X^{\roman t\eta'})=J$ for all
$X\in\im\rho$. We see that $J$ intertwines the absolutely irreducible
representations $\rho$ and $g\mapsto(\rho(g)^{-1})^{\roman t\eta'}$ and
hence that $J$ is uniquely determined up to a scalar multiple.

Since $\epsilon I\in\im\rho$ we have immediately that 
$\epsilon\epsilon^{\eta'}=1$. Replacing $K'$ by $k'(\epsilon)$, where 
$k'$ is the fixed field of $\eta'$, permits us to conclude that 
if $p=2$ then $\eta'$ is the identity (since $\epsilon=-1$ in this 
case) and if $p$ is odd then $\eta'$ has order~2 and 
inverts~$\epsilon$. It is easily checked from the formulas above that
$\rho_D$ preserves a symmetric bilinear form $f'_D$ with matrix 
$J_D=\left(\smallmatrix1&0\\0&1\endsmallmatrix\right)$ and $\rho_Q$ 
preserves an alternating bilinear form $f'_Q$ with matrix 
$J_Q=\left(\smallmatrix0&1\\-1&0\endsmallmatrix\right)$. If $p\ne2$ and 
$K'$ has an automorphism $\eta'$ inverting $\epsilon$, then $\rho_E$ 
preserves $f'_E$, the standard $\eta'$-Hermitian form (having matrix $J=I$).
Note that if $K'=\F_{q'}$ where $q'=(p')^r$, then $K'$ has a primitive
$p^{\roman{th}}$ root of~1 and an automorphism inverting it, if and only 
if $r$ is even and $p$ divides $(p')^{r/2}+1$. Given $p$ and $p'$, this 
happens if and only if $\roman{ord}_p(p')$, the least $k$ such that 
$(p')^k\equiv 1$ (mod $p$), is even.

We turn now to representations of the extraspecial groups $D^n$, 
$D^{n-1}Q$ and $E^n$. If $E_1\cdots E_n$ is a central product of 
its subgroups $E_1,\,\ldots,\,E_n$ and $\rho_1,\,\ldots,\,\rho_n$ are
absolutely irreducible representations of the $E_i$ lying over the same 
one-dimensional representation of $E_1\cap\cdots\cap E_n$, then
$\rho=\rho_1\otimes\cdots\otimes\rho_n$ (defined by 
$\rho(g)=\rho(g_1)\otimes\cdots\otimes\rho(g_n)$ whenever 
$g=g_1\cdots g_n$ with $g_i\in E_i$) is an absolutely irreducible 
representation of $E_1\cdots E_n$. Furthermore, if $\rho_i$ preserves a 
form $f_i$ then the unique bilinear form $f=f_1\otimes\cdots\otimes f_n$
which satisfies$$\postdisplaypenalty=10000
f(u_1\otimes\cdots\otimes u_n,v_1\otimes\cdots\otimes v_n)=
\prod_{i=1}^nf_i(u_i,v_i)$$is preserved by~$\rho$. Note that if each 
$f_i$ is nondegenerate then so is $f$. Thus the representation
$\bigotimes_{i=1}^n\rho_D$ of $D^n$ preserves the symmetric bilinear form
$\bigotimes_{i=1}^nf'_D$, the representation
$(\bigotimes_{i=1}^{n-1}\rho_D)\otimes\rho_Q$ of $D^{n-1}Q$ preserves the 
alternating bilinear form $(\bigotimes_{i=1}^{n-1}f'_D)\otimes 
f'_Q$, and (provided $K'$ has a suitable automorphism $\eta'$) the
representation $\bigotimes_{i=1}^n\rho_E$ of $E^n$ preserves the 
$\eta'$-Hermitian form $\bigotimes_{i=1}^nf'_E$. In each case the form is 
nondegenerate, and the degree of the representation is~$p^n$.

The faithful absolutely irreducible representation of an extraspecial
2-group over a field of odd characteristic $p'$ is unique up to equivalence, 
and we have shown that it may be realized over the prime field~$\F_{p'}$.
For odd $p$ the extraspecial $p$-group $E^n$ has an absolutely 
irreducible representation $\rho(\epsilon)$ for each choice of 
$\epsilon$ (the primitive $p^{\roman{th}}$ root of~1), and the $p-1$ 
possible choices of $\epsilon$ yield inequivalent representations. These
may all be realized over $K'=\F_{p'}(\root p\of1)=\F_{q'}$, where
$q'=(p')^r$ is the least power of $p'$ such that $p$ divides $q'-1$. A 
$p^n$-dimensional vector space over $K'$ is an $rp^n$-dimensional vector
space over~$\F_{p'}$; hence the $K'$-representation $\rho(\epsilon)$ of 
$E^n$ becomes an $\F_{p'}$-representation $\rho$ of degree~$rp^n$. The
$K'$-representation $\rho_{K'}$ obtained from $\rho$ by field extension
splits into the $r$ algebraically conjugate constituents 
$\rho(\epsilon_i)$, where $\epsilon=\epsilon_1$, $\epsilon_2$, \dots, 
$\epsilon_r$ are the algebraic conjugates of~$\epsilon$. It follows that
$\rho$ is irreducible. Note also that if $\rho$ preserves a nonzero
$\F_{p'}$-bilinear form then $\rho$ must be equivalent to its 
contragredient $\rho^*\colon g\mapsto\rho(g^{-1})^{\roman t}$ (since the
matrix of the form will intertwine $\rho$ and $\rho^*$), and since the
absolutely irreducible constituents of $\rho^*$ are the
$\rho(\epsilon_i^{-1})$ it follows that $\epsilon^{-1}$ is an algebraic
conjugate of $\epsilon$. So this can happen only when there is a field
automorphism inverting $\epsilon$, in which case, as we have seen, there
is an Hermitian form over~$K'$ preserved by~$\rho(\epsilon)$.

In the inductive step of our construction of extraspecial towers we will
embed a group $G$, which has a normal extraspecial $p$-subgroup $E(f)$, in
the automorphism group of an extraspecial $p'$-group $E(\hat f)$, by means
of a representation of $G$ which preserves $\hat f$. We prefer to use 
a representation of $G$ which is an extension of a faithful irreducible 
representation of $E(f)$. However, since $\hat f-(\hat f)^{\roman t}$ must
be nondegenerate, this is clearly impossible if $E(f)\cong D^n$ (when the 
faithful irreducible representation of $E(f)$ preserves only a symmetric
form) or if $E(f)\cong E^n$ and $\roman{ord}_p(p')$ is odd (when there is
no form at all preserved by the irreducible 
representations of $E(f)$). In these cases we are forced to resort to 
non-irreducible representations, and use the direct sum of an absolutely
irreducible representation $\tilde\rho$ of $G$ and its contragredient.
Observe that the equation
$$\pmatrix X&0\\0&X^{-\roman t}\endpmatrix\pmatrix A&B\\C&D\endpmatrix
\pmatrix X&0\\0&X^{-\roman t}\endpmatrix^{\roman t}=
\pmatrix A&B\\C&D\endpmatrix$$holds for all $X\in\im\tilde\rho$ if and
only if $A$ and $D$ are zero and $B$ and $C$ scalar multiples of~$I$.
Choosing $B=0$ and $C=I$ guarantees that the form $f'$ with matrix
$\left(\smallmatrix A&B\\C&D\endsmallmatrix\right)$ has the property 
that $f'-(f')^{\roman t}$ is nondegenerate.

\heading 4. Extraspecial groups acting on extraspecial groups\endheading

\noindent Let $V'$ be a vector space over $K'=\F_{q'}$ and
$\rho\colon E(f)\to\gl(V')$ a representation. We wish to embed $E(f)$
in the automorphism group of another extraspecial group. For this 
purpose we require a bilinear form $V'\times V'\to\F_{p'}$, where $p'$ 
is the characteristic of $K'$, rather than an $\eta'$-sesquilinear form
$V'\times V'\to K'$.

If $q'=(p')^r$ then the trace map $T\colon x\mapsto\sum_{i=0}^{r-1}x^{(p')^i}$
is a nonzero $\F_{p'}$-linear map $K'\to\F_{p'}$, and all other such 
maps are given by $x\mapsto T(\lambda x)$ for nonzero elements 
$\lambda\in K'$. Note that if $\eta'$ is any automorphism of $K'$ then
$T(x^{\eta'})=T(x)$ for all $x\in K'$. If $f'\colon V'\times V'\to 
K'$ is $\eta'$-sesquilinear and $0\ne\lambda\in K'$ then $\hat f=\hat 
f_\lambda$ defined by$$\hat f(x,y)=T(\lambda f'(x,y))$$
is an $\F_{p'}$-bilinear form $V'\times V'\to\F_{p'}$. Since we wish to
form a group $E(\hat f)$ we will require $\hat f-(\hat f)^{\roman t}$ to
be nondegenerate. Three cases will arise, as follows:
\roster
\item $\rho$ is an absolutely irreducible representation of
$E(f)\cong D^{n-1}Q$ preserving the alternating form $f'$;
\item either $E(f)\cong D^n$ or $E(f)\cong E^n$ and $\roman{ord}_p(p')$ 
is odd, and $\rho$ is the sum of an absolutely irreducible 
representation and its contragredient;
\item $p$ is odd and $\roman{ord}_p(p')$ is even, and $\rho$
is absolutely irreducible and preserves the Hermitian form $f'$.
\endroster
In the first two of these cases $\eta'=1$, while in the third it has order~2.

We define $G(\hat f)$ to be the group of all $K'$-linear transformations
of $V'$ which preserve $\hat f$:
$$G(\hat f)=\{\,g\in\gl(V')\mid \hat f(x_1g,x_2g)=\hat f(x_1,x_2)\hbox
{ for all $x_1,x_2\in V'$}\,\}.$$
It is straightforward to prove that $G(\hat f)=G(f')$, the group of all 
$g\in\gl(V')$ which preserve~$f'$.

\noindent{\it Case 1.}\indent Suppose that $p'$ is odd and $r=1$, so 
that $K'=\F_{p'}$, and $f'$ is a nondegenerate alternating form. We
simply choose $\lambda=1$, giving $\hat f=f'$ and $E(\hat f)\cong E^n$.
The group $G(f')$ is isomorphic to $\roman{Sp}(2n,p')$.

\noindent{\it Case 2.}\indent Suppose that $V'\cong W^*\oplus W$, where 
$W^*$ is the dual space of $W$, and that $f'$ is defined by
$$f'((\alpha,x),(\beta,y))=x\beta$$for all $x,y\in W$ and 
$\alpha,\beta\in W^*$. Thus $f'$ is bilinear and has matrix
$\left(\smallmatrix0&0\\I&0\endsmallmatrix\right)$ relative to a basis
comprising a basis of $W$ and the corresponding dual basis of $W^*$.
Putting $\lambda=1$ gives
$$(\hat f-(\hat f)^{\roman t})(x,y)=T((f'-(f')^{\roman t})(x,y))$$ and 
since $f'-(f')^{\roman t}$ is nondegenerate it follows (by an argument
similar to one used in the proof of 4.2 below) that
$\hat f-(\hat f)^{\roman t}$ is also nondegenerate. It is easily checked that
$$\pmatrix A&B\\C&D\endpmatrix\pmatrix 0&0\\I&0\endpmatrix
\pmatrix A&B\\C&D\endpmatrix^{\roman t}=\pmatrix 0&0\\I&0\endpmatrix$$
if and only if $B=C=0$ and $D=(A^{-1})^{\roman t}$, so that $W$ and 
$W^*$ are both $G(f')$-invariant, and $G(f')\cong\gl(W)$.
\proclaim{4.1 Proposition}If $n=\frac12r\dim V'$ (where 
$r=[K':\F_{p'}]$) then
$$E(\hat f)\cong\cases E^n&\hbox{if $p'$ is odd,}\\D^n&\hbox{if $p'=2$.}
\endcases$$\endproclaim
\demo{Proof}In the odd case the result is immediate, since $E(\hat f)$ 
has exponent~$p'$ and the dimension of $V'$ as a vector space over 
$\F_{p'}$ is~$2n$. In the case $p'=2$ all that remains is to count the
number of $v\in V'$ such that $\hat f(v,v)=0$. Now if $v=(\alpha,x)$ 
(where $x\in W$ and $\alpha\in W^*$) then $\hat f(v,v)=T(x\alpha)$.
This is zero whenever $x=0$, and for $x\ne0$ it is zero for exactly half
of the possible values of $\alpha$. Since $|W|=|W^*|=2^n$ this gives the
total number of such $v$ as $2^n+(2^n-1)2^{n-1}$. Since this equals
$2^{n-1}(2^n+1)$ the desired conclusion follows from Proposition~2.1.\enddemo

\noindent{\it Case 3.}\indent Suppose that $\eta'$ is an automorphism 
of $K'$ of order~2, and $f'$ is $\eta'$-Hermitian and nondegenerate. 
Thus, $G(f')$ is a unitary group. We choose any $\lambda$ such that 
$\lambda-\lambda^{\eta'}\ne0$.
\proclaim{4.2 Proposition}The bilinear form
$\hat f_\lambda-(\hat f_\lambda)^{\roman t}$ is nondegenerate, and
$$E(\hat f_\lambda)\cong\cases E^{sd}&\hbox{if $p'$ is odd,}\\
(D^{s-1}Q)^d&\hbox{if $p'=2$,}\endcases$$
where $2s=r=[K':\F_{p'}]$ and $d=\dim V'$.\endproclaim
\demo{Proof}For $x,y\in V'$ we have
$$\align(\hat f-(\hat f)^{\roman t})(x,y)&=T(\lambda f'(x,y))-T(\lambda 
f'(y,x))\\&=T(\lambda f'(x,y))-T(\lambda^{\eta'} f'(y,x)^{\eta'})\\
&=T(\lambda f'(x,y))-T(\lambda^{\eta'} f'(x,y))\\
&=T((\lambda-\lambda^{\eta'})f'(x,y)).\endalign$$Suppose that $y$ is such
that this is zero for all~$x$. The image of the map
$x\mapsto(\lambda-\lambda^{\eta'})f(x,y)$ is a $K'$-subspace of $K'$, and
not equal to $K'$ since it is contained in the kernel of~$T$. So it must
be zero, and we conclude that $f'(x,y)=0$ for all~$x$. Since $f'$ is 
nondegenerate, $y=0$. Hence $\hat f_\lambda-(\hat f_\lambda)^{\roman t}$
is nondegenerate.

We need the following simple proposition, which will also be used in
Section~6.
\proclaim{4.3 Proposition}Let $k'=\{\,a\in K'\mid a^{\eta'}=a\,\}$. Then 
there is a nonzero $t_0\in K'$ such that
$k't_0=C=\{\,a\in K'\mid a^{\eta'}=-a\,\}$. Furthermore, if $0\ne t\in K'$ 
then $k't$ is contained in $\ker T$ if and only if $k't=C$.\endproclaim
\demo{Proof}Observe that $a\mapsto a^{\eta'}+a$ is a nonzero $k'$-linear map
from $K'$ to $k'$. Since the dimension of $K'$ over $k'$ is 2, its kernel
$C$ is a 1-dimensional subspace, and hence equals $k't_0$ for some
$t_0\in K'$. Clearly $C\subseteq\ker T$. If $k't\ne k't_0$ then
$K'=k't_0\oplus k't$, which precludes $k't\subseteq\ker T$ since $T$ is not
the zero map.\enddemo
It is easily proved (cf. \cite{\hup, p. 235}) 
that $V'$ has a $K'$-basis which is orthonormal relative to
$f'$, so that $f'$ can be written as the direct sum of $d$ copies of
the form $(x,y)\mapsto xy^{\eta'}$ on the one-dimensional space~$K'$.
Correspondingly, $\hat f_\lambda$ is the direct sum of $d$ copies of
$\hat f_0\colon(x,y)\mapsto T(\lambda xy^{\eta'})$. Thus
$E(\hat f_\lambda)\cong E(\hat f_0)^d$.

If $p'$ is odd then $E(\hat f_0)\cong E^s$, since the dimension of $K'$
over $\F_{p'}$ is $2s$. When $p'=2$ our task is to count the number of 
$x\in K'$ such that $T(\lambda xx^{\eta'})=0$. Since 
$\lambda\ne\lambda^{\eta'}=-\lambda^{\eta'}$, Proposition~4.3 shows that
$k'\lambda\not\subseteq\ker T$. Note that $xx^{\eta'}\in k'$ 
for all $x\in K'$. Now $xx^{\eta'}=0$ if and only if $x=0$, while 
$xx^{\eta'}=y$ has $2^s+1$ solutions whenever $0\ne y\in k'$. Since 
$y\mapsto T(\lambda y)$ is a nonzero $\F_2$-linear map $k'\to\F_2$, its 
kernel must have $2^{s-1}-1$ nonzero elements. Hence there are 
$(2^{s-1}-1)(2^s+1)$ nonzero $x\in K'$ such that $T(\lambda xx^{\eta'})=0$.
Therefore the total number of solutions of $T(\lambda xx^{\eta'})=0$ is
$$1+(2^{s-1}-1)(2^s+1)=2^{s-1}(2^s-1)$$and by Proposition 2.1 it follows
that $E(\hat f_0)\cong D^{s-1}Q$.\enddemo

Note that since $DD=QQ$ it follows that
$$(D^{s-1}Q)^d\cong\cases D^{ds}&\hbox{if $d$ is 
even,}\\D^{ds-1}Q&\hbox{if $d$ is odd,}\endcases$$
completing the proof of 4.2.
\heading 5. Constructing the factor set\endheading

\noindent Let $V$ be a vector space over a field $K$ of characteristic $p$,
and let $f\colon V\times V\to\F_p$ be a $\F_p$-bilinear map such that
$f-f^{\roman t}$ is nondegenerate. (The ``$f$'' of this section 
corresponds to the ``$\hat f$'' of the previous section.) Let the 
dimension of $V$ over $\F_p$ be $2n$.

Let $V'$ be a $p^n$-dimensional vector space over a field $K'$ of
characteristic~$p'$, and let $\rho\colon E(f)\to\gl(V')$ 
be an absolutely irreducible representation. In this section we will 
extend $\rho$ to a projective representation of $G(f)\semiprod E(f)$ and
calculate the factor set involved. In the next section we 
will show that the factor set splits (at least in the cases that concern
us), so that $\rho$ extends to a representation of $G(f)\semiprod E(f)$.
Furthermore, we will show that if $\rho$ preserves a form $f'$ on $V'$ then
the extension of $\rho$ also preserves $f'$.

Observe that $x\mapsto\rho(x,0)$ is a projective representation of the
additive group of~$V$, since for all $x,y\in V$,
$$\align\rho(x,0)\rho(y,0)&=\rho(x+y,f(x,y))\\
&=\rho(0,f(x,y))\rho(x+y,0)\\
&=\epsilon^{f(x,y)}\rho(x+y,0)\endalign$$
where $\epsilon\in K'$ is a primitive $p^{\roman{th}}$ root of~1. 
Since $\rho$ is absolutely irreducible, so too is this projective 
representation.

To minimize the use of superscripts, we define $\exp(t)=\epsilon^t$.
\proclaim{5.1 Lemma}Let $G$ be any group and let $\rho_i\colon G\to\gl(V_i)$
(for $i=1,2$) be projective representations with the same factor set 
$\alpha$. If $a\colon V_1\to V_2$ is an arbitrary linear map, then 
$$s=\sum_{y\in G}\rho_1(y)a\rho_2(y)^{-1}$$satisfies
$$\rho_1(x)s=s\rho_2(x)\qquad\hbox{for all $x\in G$.}$$\endproclaim
\demo{Proof}We have
$$\align \rho_1(x)s&=\sum_{y\in G}\alpha(x,y)\rho_1(xy)a\rho_2(y)^{-1}
\quad\hbox{(as $\rho_1(x)\rho_1(y)=\alpha(x,y)\rho_1(xy)$)}\\
&=\sum_{y\in G}\rho_1(xy)a\rho_2(xy)^{-1}\rho_2(x)\\
\noalign{\vskip-\smallskipamount\rightline{(as
$\rho_2(xy)^{-1}\rho_2(x)=\alpha(x,y)\rho_2(y)^{-1}$)}}
&=s\rho_2(x).\endalign$$\enddemo
Observe that if $g\in G(f)$ then $x\mapsto\rho(xg,0)$ is also a 
projective representation of the additive group of $V$, and it has the 
same factor set as $x\mapsto\rho(x,0)$. Following Ward [\wardi] we define
$$s(g)=|V|^{-1}\sum_{y\in V}\rho(y,0)\rho(yg,0)^{-1},$$
so that (by Lemma 5.1) $s(g)$ intertwines these projective 
representations. Since they are irreducible it follows from Schur's 
Lemma that $s(g)$ is either zero or invertible.
As $\rho(yg,0)^{-1}=\exp(f(y,y))\rho(-yg,0)$ we find that
$$s(g)=|V|^{-1}\sum_{y\in V}\exp(f(y,y(1-g)))\rho(y(1-g),0).\tag\$$$

If $F$ is a bilinear form on $V$ then given a subspace $W$ of $V$ we 
define
$$\eqalignno{W^F&=\{\,y\in V\mid F(x,y)=0\hbox{ for all $x\in W$}\,\}\cr
\noalign{\hbox{and}}
^FW&=\{\,x\in V\mid F(x,y)=0\hbox{ for all $y\in W$}\,\}.\cr}$$
For each $g\in\gl(V)$ we define $K(g)$ and $I(g)$ to be (respectively)
the kernel and image of $1-g$.
\proclaim{5.2 Lemma}Suppose that $g\in\gl(V)$ preserves the bilinear 
form~$F$. Then $K(g)$ is contained in both $^FI(g)$ and $I(g)^F$. 
Furthermore, if $F$ is nondegenerate then 
$^FI(g)=I(g)^F=K(g)$.\endproclaim
\demo{Proof}Let $x\in K(g)$. If $v\in V$ then
$$\align F(v(1-g),x)&=F(v,x)-F(vg,x)\\
&=F(vg,xg)-F(vg,x)\qquad\hbox{(since $g$ preserves $F$)}\\
&=F(vg,xg-x)\\
&=0\qquad\qquad\qquad\hbox{(since $xg-x=-x(1-g)=0$).}\endalign$$
Thus $F(y,x)=0$ for all $y\in I(g)$, so that $x\in I(g)^F$. So 
$K(g)\subseteq I(g)^F$. If $F$ is nondegenerate then
$$\dim I(g)^F=\dim V-\dim I(g)=\dim K(g),$$
and we deduce that $I(g)^F=K(g)$.

The corresponding facts concerning $^FI(g)$ can be proved by similar 
arguments.\enddemo

Let $x,y\in I(g)$. If $u,u'\in V$ are such that
$x=u(1-g)=u'(1-g)$ then $u'=u+v$ with $v\in K(g)$, and by 5.2 we have
$$f(u',y)=f(u,y)+f(v,y)=f(u,y).$$Hence, following Wall [\wall], we may
define $f_g$ on $I(g)\times I(g)$ by
$$f_g(x,y)=f(u,y)\qquad\hbox{for all $u$ such that $x=u(1-g)$.}\tag\dag$$
Clearly $f_g$ is $\F_p$-bilinear. Note, furthermore, that if $y=v(1-g)$ then
$$\aligned f_g(x,y)&=f(u,v-vg)=f(u,v)-f(u,vg)\\
		 &=f(ug,vg)-f(u,vg)=-f(x,vg)=f(x,y-v).
\endaligned\tag\ddag$$

If $y\in V$ then $f(y,y(1-g))=f_g(y(1-g),y(1-g))$, and we can rewrite 
the formula $(\$)$ above as
$$s(g)=|I(g)|^{-1}\sum_{x\in I(g)}\exp(f_g(x,x))\rho(x,0).$$
\proclaim{5.3 Proposition}For each $g\in G(f)$ the transformation $s(g)$
is invertible.\endproclaim
\demo{Proof}Since $\rho$ is absolutely irreducible its enveloping 
algebra (the linear span of the set
$\{\,\rho(x,z)\mid x\in V,\,z\in\F_p\,\}$) is the $p^{2n}$-dimensional 
space of all linear transformations of $V'$. As $\rho(x,z)$ is a scalar
multiple of $\rho(x,0)$ it follows that the $p^{2n}$ transformations 
$\rho(x,0)$ must be linearly independent. It now follows immediately 
from our expression for $s(g)$ that $s(g)\ne0$, and therefore (by 
Schur's Lemma) that $s(g)$ is invertible.\enddemo
\proclaim{5.4 Theorem}With the notation as above
$$g(x,z)\mapsto s(g)\rho(x,z)$$
defines a projective representation of $G(f)\semiprod E(f)$. Furthermore, if 
$g,\,h\in G(f)$ then $s(g)s(h)=\sigma(g,h)s(gh)$, where
$$\sigma(g,h)=|I(g)|^{-1}|I(h)|^{-1}|I(gh)|\sum_{x\in I(g)\cap I(h^{-1})}
\exp(\gamma_{g,h}(x,x))$$the $\F_p$-bilinear form $\gamma_{g,h}$ on
$I(g)\cap I(h^{-1})$ being given by
$$\gamma_{g,h}(u,v)=f_g(u,v)-f_{h^{-1}}(u,v).$$\endproclaim
\demo{Proof}For each $g\in G(f)$ we have 
$s(g)^{-1}\rho(x,z)s(g)=\rho(xg,z)$ for all $x\in V$ and $z\in\F_p$, and 
it follows that if $g,h\in G(f)$ then
$$s(h)^{-1}s(g)^{-1}\rho(x,z)s(g)s(h)=s(gh)^{-1}\rho(x,z)s(gh)$$
for all $x$ and $z$. By Schur's Lemma $s(g)s(h)=\sigma(g,h)s(gh)$ for 
some scalar $\sigma(g,h)$.

Multiplying the expressions for $s(g)$ and $s(h)$ we find that 
$s(g)s(h)$ is the sum of all terms
$$|I(g)|^{-1}|I(h)|^{-1}\exp(f_g(x,x)+f_h(y,y)+f(x,y))\rho(x+y,0)$$
for $x\in I(g)$ and $y\in I(h)$. The coefficient of $\rho(0,0)$ in 
$s(g)s(h)$ is therefore
$$|I(g)|^{-1}|I(h)|^{-1}\sum_{x\in I(g)\cap I(h)}\exp(f_g(x,x)+
f_h(-x,-x)+f(x,-x)).$$If $x=u(1-h)$ then $x=(-uh)(1-h^{-1})$. Hence
$I(h)=I(h^{-1})$, and furthermore,
$$f_{h^{-1}}(x,x)=f(-uh,x)=f(x-u,x)=f(x,x)-f_h(x,x).$$Thus the 
coefficient obtained above can be written as
$$|I(g)|^{-1}|I(h)|^{-1}\sum_{x\in I(g)\cap I(h^{-1})}\exp(\gamma_{g,h}(x,x)).$$
Since this must equal the coefficient of $\rho(0,0)$ in $\sigma(g,h)s(gh)$,
which is $\sigma(g,h)|I(gh)|^{-1}$, the result follows.( Note the
similarity with (4.3.1) of \cite{\brw}.)\enddemo
To show that $\rho$ extends to a representation of $G(f)\semiprod E(f)$, 
rather than merely a projective representation, it is necessary to find 
for each $g\in G(f)$ a nonzero $\mu(g)\in K'$ such that
$\sigma(g,h)=\mu(g)\mu(h)\mu(gh)^{-1}$, so that defining 
$s'(g)=\mu(g)^{-1}s(g)$ leads to $s'(g)s'(h)=s'(gh)$.
Our next two theorems will be used later in the proof that such a 
function~$\mu$ exists.
\proclaim{5.6 Theorem}The radical of the quadratic form 
$x\mapsto\gamma_{g,h}(x,x)$ on $I(g)\cap I(h^{-1})$ is
$$R_{g,h}=\{\,x\in V\mid x=u(1-g)=u(1-h^{-1})\hbox{ for some $u\in V$}\,\}.$$
\endproclaim
\demo{Proof}Let $\gamma=\gamma_{g,h}$. The radical of the quadratic form 
consists of those $x$ such that $\gamma(x+y,x+y)=\gamma(y,y)$ for all~$y$,
or, equivalently, those $x$ such that $\gamma(x,x)=0$ and
$(\gamma+\gamma^{\roman t})(x,y)=0$ for all~$y$.

If $x\in R_{g,h}$ (defined as above) then 
$f_g(x,y)=f_{h^{-1}}(x,y)$ (by (\dag)) and $f_g(y,x)=f_{h^{-1}}(y,x)$
(by (\ddag)) for all $y\in I(g)\cap I(h^{-1})$, and it
follows that $x$ is in the radical.
Conversely, let $x=u(1-g)=w(1-h^{-1})$ be an element of the radical. Let
$y$ be an arbitrary element of $I(g)\cap I(h^{-1})$. By (\dag) we have
$\gamma(x,y)=f(u,y)-f(w,y)=f(u-w,y)$, while by (\ddag)
$$\align\gamma^{\roman t}(x,y)&=\gamma(y,x)\\
&=f(y,x-u)-f(y,x-w)\\
&=f(y,w-u).\endalign$$
Thus $(f-f^{\roman t})(u-w,y)=(\gamma+\gamma^{\roman t})(x,y)=0$, and 
this holds for all $y\in I(g)\cap I(h^{-1})$. Since $F=f-f^{\roman t}$ is 
nondegenerate and preserved by $g$ and $h$ it follows from 5.2 that
$$(I(g)\cap I(h^{-1}))^F=I(g)^F+I(h^{-1})^F=K(g)+K(h^{-1}).$$
Hence $u-w=u_0-w_0$ for some $u_0\in K(g)$ and $w_0\in K(h^{-1})$.
Putting $u'=u-u_0=w-w_0$ we see that $x=u'(1-g)=u'(1-h^{-1})$. Hence 
$x\in R_{g,h}$.\enddemo
Theorem 5.6 will be used to compute an expression for $\sigma(g,h)$, and
the following theorem will be used to split~$\sigma$.
\proclaim{5.7 Theorem}If $g,h\in G(f)$ then 
$$\dim_K(I(g)\cap I(h^{-1}))+\dim_KR_{g,h}=i(g)+i(h)-i(gh)$$
where we have defined $i(k)=\dim_KI(k)$ for all $k\in G(f)$.\endproclaim
\demo{Proof}Note that elements of $G(f)$ are $K$-linear transformations,
so that all the spaces involved in the theorem statement are indeed
$K$-subspaces of~$V$. It suffices, however, to prove the corresponding 
statement for $\F_p$-dimensions, since the desired conclusion will then 
follow by dividing by $[K:\F_p]$.

Recall that $\dim V=2n$. Since $F=f-f^{\roman t}$ is nondegenerate we 
have
$$\align\dim(I(g)&\cap I(h^{-1}))=2n-\dim(I(g)\cap I(h^{-1}))^F\\
&=2n-\dim(I(g)^F+I(h^{-1})^F)\\
&=2n-\dim K(g)-\dim K(h^{-1})+\dim(K(g)\cap K(h^{-1}))\\
&=i(g)-(2n-i(h^{-1}))+\dim(K(g)\cap K(h^{-1}))\endalign$$
by Lemma 5.2. By Theorem 5.6 we have $R_{g,h}=K(gh)(1-g)$; however, it is 
easily shown that $\ker(1-g)\cap K(gh)=K(g)\cap K(h^{-1})$, and therefore
$$\align\dim R_{g,h}&=\dim K(gh)-\dim(K(g)\cap K(h^{-1}))\\
&=2n-i(gh)-\dim(K(g)\cap K(h^{-1})).\endalign$$
Clearly $i(h^{-1})=i(h)$; so adding our formulas for
$\dim(I(g)\cap I(h^{-1}))$ and $\dim R_{g,h}$ gives the required result.
\enddemo

Suppose now that $\rho$ preserves a nonzero $\eta'$-sesquilinear form, 
$\eta'$ being an automorphism of $K'$ which inverts~$\epsilon$. For 
simplicity we regard $\rho$ as a matrix representation, and we let $J$ 
be the matrix of the form. Then $J$ is nonsingular (since $\rho$ is 
absolutely irreducible) and$$J^{-1}\rho(y,0)J=
(\rho(y,0)^{\roman t\eta'})^{-1}$$
for all $y\in V$. From the definition of $s(g)$ it follows that
$$\align J^{-1}s(g)J&=|V|^{-1}\sum_{y\in V}(J^{-1}\rho(y,0)J)
(J^{-1}\rho(yg,0)^{-1}J)\\
&=|V|^{-1}\sum_{y\in V}(\rho(y,0)^{-1})^{\roman t\eta'}
\rho(yg,0)^{\roman t\eta'}\\
&=|V|^{-1}\sum_{y\in V}(\rho(yg,0)\rho(y,0)^{-1})^{\roman t\eta'}\\
&=\bigr(|V|^{-1}\sum_{x\in V}\rho(x,0)
\rho(xg^{-1},0)^{-1}\bigl)^{\roman t\eta'}\\
&=s(g^{-1})^{\roman t\eta'}.\endalign$$
Hence the following theorem holds.
\proclaim{5.8 Theorem}If there exists a function $\mu\colon G(f)\to K'$ 
satisfying 
$$\sigma(g,h)=\mu(g)\mu(h)\mu(gh)^{-1}$$ 
and 
$\mu(g)=\mu(g^{-1})^{\eta'}$ for all $g,h\in G(f)$ then 
$s'(g)=\mu(g)^{-1}s(g)$ defines an extension of $\rho$ which preserves 
any $\eta'$-sesquilinear form preserved by~$\rho$.\endproclaim

\heading 6. Splitting the factor set\endheading

\noindent Let $f$ and $V$ be as in the previous section; our aim is to find a 
function $\mu$ such that $\sigma(g,h)=\mu(g)\mu(h)\mu(gh)^{-1}$ and 
$\mu(g)=\mu(g^{-1})^{\eta'}$ for all $g,h\in G(f)$. We treat
three separate cases, corresponding to the three cases in Section~4: 
\roster
\item $p$ is odd, $K=\F_p$ and $f^{\roman t}=-f$;
\item $V=W^*\oplus W$ where $W$ is a vector space over $K$ of 
$\F_p$-dimension~$n$, and $f$ is defined by
$$\quad f((\alpha,x),(\beta,y))=T(x\beta)\quad\hbox{for all $x,y\in W$
and $\alpha,\beta\in W^*$}$$where $T\colon K\to\F_p$ is the trace map;
\item $K$ has an automorphism $\eta$ of order~2 and
$f=f_\lambda$ is given by
$$f(u,v)=T(\lambda F(u,v))\quad\hbox{for all $u,v\in V$}$$
where $F\colon V\times V\to K$ is a $\eta$-Hermitian form, $T\colon K\to\F_p$
is the trace map, and $\lambda$ is some fixed element of $K$ with 
$\lambda\ne\lambda^\eta$.
\endroster
Investigating the second of these cases first, assume the hypotheses (2)
above, and let $|K|=q$. We identify $W^*$ and $W$ with subspaces $W_1$ and
$W_2$ of $V$ in the obvious fashion, so that elements of $V$ have the form
$v_1+v_2$ with $v_1\in W_1=W^*$ and $v_2\in W_2=W$. We showed in Section~4
that $W$ and $W^*$ are both 
$G(f)$-invariant, and it follows readily that
$$I(g)=I_1(g)\oplus I_2(g)\quad\hbox{ for all $g\in G(f)$}$$
where $I_1(g)=I(g)\cap W_1$ and $I_2(g)=I(g)\cap W_2$. Furthermore, if 
$g,h\in G(f)$ then $$I(g)\cap I(h^{-1})=
(I_1(g)\cap I_1(h^{-1}))\oplus(I_2(g)\cap I_2(h^{-1})),$$
and it follows that
$$\sum_{x\in I(g)\cap I(h^{-1})}\exp(\gamma_{g,h}(x,x))=
\sum_{x_1\in I_1}\sum_{x_2\in I_2}\exp(\gamma_{g,h}(x_1+x_2,x_1+x_2))$$
where $\gamma_{g,h}$ is as defined in Theorem~5.4 and $I_i=I_i(g)\cap 
I_i(h^{-1})$ for each~$i$. (In view of the formula
for $\sigma(g,h)$ in Theorem 5.4, our task is to evaluate this sum.)
\proclaim{6.1 Theorem}Assume the conditions of the above preamble, and 
define $j(g)=\dim(I(g)\cap W)$ for all $g\in G(f)$. Then the function $\mu$
defined on $G(f)$ by $\mu(g)=q^{-j(g)}$ splits the factor set $\sigma$ and
satisfies $\mu(g)=\mu(g^{-1})^{\eta'}$ for all~$g\in G(f)$. \endproclaim
\demo{Proof} Since the map $(x,\beta)\mapsto x\beta$ from
$W\times W^*$ to $K$ is nondegenerate and $G(f)$-invariant, reasoning
parallel to that used to prove Lemma 5.2 shows that for all $g\in G(f)$,
$$(W^*(1-g))^\perp=W\cap\ker(1-g)$$
(where $(W^*(1-g))^\perp$ is defined as the set of all $v\in W$ which are
annihilated by all $\beta\in W^*(1-g)=I_1(g)$). Observe that, as a consequence, 
$I_1(g)$ and $I_2(g)$ have the same dimension. Note also that if $g$ and
$h$ are both in $G(f)$ then$$(I_1(g)\cap I_1(h^{-1}))^\perp=
(W_2\cap\ker(1-g))+(W_2\cap\ker(1-h^{-1})).$$

Let $x_i\in I_i(g)\cap I_i(h^{-1})$ (for $i=1,2$), and let
$u_1,v_1\in W_1$ and $u_2,v_2\in W_2$ with $x_i=u_i(1-g)=v_i(1-h^{-1})$ 
for each $i$. We find that
$$\align\gamma_{g,h}(x_1+x_2,x_1+x_2)&
=f_g(x_1+x_2,x_1+x_2)-f_{h^{-1}}(x_1+x_2,x_1+x_2)\\
&=f(u_1+u_2,x_1+x_2)-f(v_1+v_2,x_1+x_2)\\
&=T((u_2-v_2)x_1).\endalign$$

If $u_2-v_2\notin (I_1(g)\cap I_1(h^{-1}))^\perp$ then 
$T((u_2-v_2)(x_1))$ assumes all values equally often as $x_1$ varies over
all elements of $I_1(g)\cap I_1(h^{-1})$, and since
$\sum_{i=0}^{p-1}\epsilon^i=0$ it follows that
$$\sum_{x_1\in I_1(g)\cap I_1(h^{-1})}\exp(T((u_2-v_2)(x_1)))=0.$$
If $u_2-v_2\in(I_1(g)\cap I_1(h^{-1}))^\perp$ then clearly this
sum is  $|I_1(g)\cap I_1(h^{-1})|$.

Now $u_2-v_2\in (I_1(g)\cap I_1(h^{-1}))^\perp$ if and only if
$u_2-v_2=u'+v'$ with $u'\in W_2\cap\ker(1-g)$ and $v'\in W_2\cap\ker(1-h^{-1})$,
and (as in the proof of Theorem 5.6) this happens if and only if
$x_2=z(1-g)=z(1-h^{-1})$ for some $z\in W_1$. Furthermore, if $R_2$ is 
the set of all such $x_2$ then we find, as in Theorem 5.7, that
$$\dim(I_1(g)\cap I_1(h^{-1}))+\dim R_2=j(g)+j(h)-j(gh).$$Thus
$$\align\sum_{x_2}\sum_{x_1}\gamma_{g,h}(x_1+x_2,x_1+x_2)&=
|I_1(g)\cap I_1(h^{-1})|\,|R_2|\\
&=q^{j(g)}q^{j(h)}q^{-j(gh)}\endalign$$
so that Theorem 5.4 gives
$$\sigma(g,h)=|I(g)|^{-1}q^{j(g)}|I(h)|^{-1}q^{j(h)}|I(gh)|q^{-j(gh)}.$$
Since $\dim I(g)=2j(g)$ (and likewise for $h$ and $gh$) it follows that 
$\mu$ splits $\sigma$. The other assertion is trivial, since 
$j(g)=j(g^{-1})$ and $q$ is fixed by all automorphisms of~$K'$.\enddemo

We investigate the third case next. As in the case just considered the 
main task is to evaluate the sum $\sum\gamma_{g,h}(x,x)$ over $x\in 
I(g)\cap I(h^{-1})$. Note that since $K$ has a nontrivial involutory 
automorphism its order is a square: $|K|=q^2$, where $q$ is the order
of the fixed field of~$\eta$.

Reasoning as in Section 5, but using the $\eta$-Hermitian form $F$ in place of
the bilinear form $f$, we see that there is a well defined $\eta$-sesquilinear
form $F_{g,h}$ on $I(g)\cap I(h^{-1})$ such that
$$F_{g,h}(x,y)=F(u-v,y)=F(x,-u'g+v'h^{-1})=F(x,v'-u')$$
whenever $x=u(1-g)=v(1-h^{-1})$ and $y=u'(1-g)=v'(1-h^{-1})$. Since $F$ 
is Hermitian we see that
$$F_{g,h}(y,x)=F(y,v-u)=-(F(u-v,y))^\eta=-(F_{g,h}(x,y))^\eta$$
so that $F_{g,h}$ is skew-Hermitian. Moreover, an argument similar to 
that used in the proof of Theorem~5.6 shows that $R_{g,h}$ (as defined 
in~5.6) is the radical of~$F_{g,h}$.
\proclaim{6.2 Lemma}Let $U$ be a $d$-dimensional vector space over the 
field~$K$, and $(\ ,\ )$ a nondegenerate skew-Hermitian form on $U$. 
Then there exists a basis $u_1,\,u_2,\,\ldots,\,u_d$ of $U$ such that
$$(u_i,u_j)=\cases 0&\hbox{if $i\ne j$}\\ z_i&\hbox{if $i=j$}\endcases$$
where the elements $z_i\in K$ satisfy $z_i^\eta=-z_i\ne0$.\endproclaim
\demo{Proof} See \cite{\hup, p. 235}.
\enddemo
Define $k=\{\,a\in K\mid a^\eta=a\,\}$, and recall that $\lambda$ is
an element of $K$ such that $\lambda\notin k$. If $z$ satisfies
$z^\eta=-z\ne0$ then $(\lambda z)^\eta\ne-\lambda z$, and
by Proposition 4.3 it follows that $k\lambda z\not\subseteq\ker T$. Now
if $0\ne a\in K$ then $0\ne aa^\eta\in k$, each nonzero element of $k$ 
occurring exactly $q+1$ times. Hence
$$\align\sum_{a\in K}\exp(T(aa^\eta\lambda z))&
=1+(q+1)\sum_{0\ne a\in k}\exp(T(a\lambda z))\\
&=-q+(q+1)\sum_{a\in k}\exp(T(a\lambda z))\\
&=-q\qquad\qquad\hbox{(since $\sum_{i=0}^{p-1}\epsilon^i=0$).}\endalign$$
Let $d=\dim(I(g)\cap I(h^{-1})/R_{g,h})$. By Lemma~6.2 we can choose 
elements $u_1,\,u_2,\,\ldots,\,u_d\in I(g)\cap I(h^{-1})$ such that
$$F_{g,h}(u_i,u_j)=\cases0&(i\ne j)\\z_i&(i=j)\endcases$$and each element 
of $I(g)\cap I(h^{-1})$ is uniquely expressible in the form 
$(\sum_{i=1}^da_iu_i)+u$ with $a_i\in K$ and $u\in R_{g,h}$. This gives
$$\displaylines{\sum_{x\in I(g)\cap I(h^{-1})}\exp(T(\gamma_{g,h}(x,x)))
=\sum_{x\in I(g)\cap I(h^{-1})}\exp(T(\lambda F_{g,h}(x,x)))\cr
=\sum_{u\in R_{g,h}}\sum_{a_1\in K}\cdots\sum_{a_d\in K}
\exp(T(\lambda(a_1a_1^\eta z_1+a_2a^\eta_2z_2+\cdots+a_da^\eta_dz_d)))\cr
=|R_{g,h}|\prod_{i=1}^d\Bigl(\sum_{a\in K}
\exp(T(aa^\eta\lambda z_i))\Bigr).\cr}$$
This in turn is equal to
$$\align (q^2)^{\dim R_{g,h}}(-q)^d&=(-q)^{2\dim R_{g,h}+d}\\
&=(-q)^{i(g)+i(h)-i(gh)}\qquad\hbox{(where $i(g)=\dim I(g)$)}\endalign$$
by 5.7, since $2\dim R_{g,h}+d=\dim R_{g,h}+\dim(I(g)\cap I(h^{-1}))$.

The above calculations and Theorem~5.4 immediately yield the following
theorem:
\proclaim{6.3 Theorem}In the situation above the function $\mu$ defined 
on $G(f)$ by $\mu(g)=(-q)^{i(g)}$ splits the factor set $\sigma$ and 
satisfies $\mu(g)=\mu(g^{-1})^{\eta'}$ for all~$g\in G(f)$.\endproclaim

We have now dealt with the third of our three cases, so that only the 
first is left. Assume, therefore, that
$p$ is an odd prime, $K=\F_p$ and $f$ is a nondegenerate alternating 
form on~$V$. For each $g\in G(f)$ we define $\delta_g=\chi(\det M)$,
where $M$ is the matrix of the form $f_g$ (relative to any basis of 
$I(g)$) and $\chi(t)$ is 1 if $t$ is a square, $-1$ otherwise. It is 
easily seen that $\delta_g$ is well-defined as changing basis multiplies
the determinant of the matrix of a form by a nonzero square.

Let $g,\,h\in G(f)$, and assume first of all that $1-gh$ is
invertible. Since $K(g)\cap K(h^{-1})\subseteq K(gh)=\{0\}$ we may
choose a basis $v_1,\,v_2,\,\ldots,\,v_{2n}$ of $V$ such that 
$v_1,\,v_2,\,\ldots,\,v_r$ is a basis of $K(g)$ and 
$v_{s+1},\,v_{s+2},\,\ldots,v_{2n}$ is a basis of $K(h^{-1})$, for some 
$r,\,s$ with $0\le r\le s\le 2n$. We identify endomorphisms of $V$ with 
their matrices relative to this basis. The division of the basis into three
parts (the first $r$ terms, the next $s-r$, and the remaining $2n-s$)
results in a corresponding partitioning of the matrices. Let $j$ be the
matrix of~$f$.

Since the first part of the basis is in the kernel of $1-g$ we see
that $1-g$ has the form
$$\pmatrix 0&0&0\\*&*&*\\*&*&*\endpmatrix$$
where the $*$'s indicate entries which are not yet relevant. Since
$$(1-g)j=j(1-g^{-\roman t})=\pmatrix 0&*&*\\0&*&*\\0&*&*\endpmatrix$$
we obtain $1-g=\left(\smallmatrix 0&0&0\\0&a&b\\0&c&d\endsmallmatrix
\right)j^{-1}$. Similarly
$1-h^{-1}=\left(\smallmatrix e&f&0\\g&h&0\\0&0&0\endsmallmatrix\right)j^{-1}$,
and therefore
$$h^{-1}-g=\pmatrix -e&-f&0\\-g&a-h&b\\0&c&d\endpmatrix j^{-1}.$$

Note that $I(g)$ consists of all vectors of the form $(0,*,*)j^{-1}$ and
$I(h^{-1})$ of all vectors of the form $(*,*,0)j^{-1}$. We have
$$\align f_g\Bigl((0,x,y)j^{-1},(0,z,w)j^{-1}\Bigr)&=
f\Bigl((0,x,y)\pmatrix 1&0&0\\0&a&b\\0&c&d\endpmatrix^{-1}\!\!\!,\,\,
(0,z,w)j^{-1}\Bigr)\\
&=(0,x,y)\!\pmatrix 1&0&0\\0&a&b\\0&c&d\endpmatrix^{-1}\!\!\!j
\Bigl(-j^{-1}\pmatrix0\\z^{\roman t}\\w^{\roman t}\endpmatrix\Bigr).\endalign$$
In effect, the matrix of $f_g$ is
$-\left(\smallmatrix a&b\\c&d\endsmallmatrix\right)^{-1}
=\left(\smallmatrix a'&b'\\c'&d'\endsmallmatrix\right)$ (say),
and $\delta_g$ is 1 if the determinant of this is a square, $-1$ if 
not. Furthermore, $a'$ is the matrix of the restriction of $f_g$ to
$I(g)\cap I(h^{-1})$.

Likewise, $f_{h^{-1}}$ has matrix
$-\left(\smallmatrix e&f\\g&h\endsmallmatrix\right)^{-1}
=\left(\smallmatrix e'&f'\\g'&h'\endsmallmatrix\right)$, and its 
restriction to $I(g)\cap I(h^{-1})$ has matrix~$h'$. Thus
$a'-h'$ is the matrix of~$\gamma_{g,h}$.

The form $\gamma_{g,h}$ is symmetric (and therefore the matrix $a'-h'$ 
is symmetric) since if $x=u(1-g)=v(1-h^{-1})$ then $(\ddag)$ and $(\dag)$ yield
$$\gamma_{g,h}(y,x)=f(y,(x-u)-(x-v))=f(y,v-u)=f(u-v,y)=\gamma_{g,h}(x,y).$$
Since the characteristic of $K$ is odd the radical of the bilinear form
$\gamma_{g,h}$ is the same as the radical of the quadratic form
$x\mapsto\gamma_{g,h}(x,x)$, which is $(\ker(h^{-1}-g))(1-g)=\{0\}$ (in 
view of our assumption that $1-gh$ is invertible). It is easily shown 
(cf. Lemma 6.2) that there is a basis for
$I(g)\cap I(h^{-1})$ relative to which the matrix of $\gamma_{g,h}$ is
diagonal, with (nonzero) entries $z_1,\,z_2,\,\dots,\,z_d$ say. It follows that
$$\sum_{x\in I(g)\cap I(h^{-1})}\exp\gamma_{g,h}(x,x)=
\prod_{i=1}^d\sum_{a\in\F_p}\exp(z_ia^2).$$
The following result is standard (cf. \cite{\wardii, Lemma 2.1}).
\proclaim{6.4 Proposition}For $0\ne z\in\F_p$ define
$\theta(z)=\sum_{a\in\F_p}\exp(za^2)$, and let $\theta=\theta(1)$.
Then $\theta(z)=\chi(z)\theta$, and $\theta^2=\chi(-1)p$.
Furthermore, if $K'$ has an automorphism $\eta'$ inverting $\epsilon$ 
(the primitive $p^{\fam0 th}$ root of~1) then 
$\theta^{\eta'}=\chi(-1)\theta$.\endproclaim
It follows from 6.4 that
$$\sum_{x\in I(g)\cap I(h^{-1})}\exp\gamma_{g,h}(x,x)
=\theta^d\chi(\det(a'-h')).$$
where $d=\dim(I(g)\cap I(h^{-1}))$. Since the radical of
$\gamma_{g,h}$ is zero, Theorem~5.7 yields $d=i(g)+i(h)-i(gh)$.

Since $1-gh$ is invertible the bilinear form $f_{gh}$ has matrix 
$(1-gh)^{-1}j$, and since $\det h=1$ it follows that
$$\delta_{gh}=\chi(\det (h^{-1}-g)j)=
\chi(\det\pmatrix -e&-f&0\\-g&a-h&b\\0&c&d\endpmatrix\,).$$
Observe that
$$\pmatrix 1&0&0\\0&a'&b'\\0&c'&d'\endpmatrix
\pmatrix -e&-f&0\\-g&a-h&b\\0&c&d\endpmatrix
\pmatrix e'&f'&0\\ g'&h'&0\\0&0&1\endpmatrix
=\pmatrix 1&0&0\\*&a'-h'&0\\*&*&-1\endpmatrix$$
so that $\chi(\det(a'-h'))=
\delta_g\delta_{h^{-1}}\delta_{gh}\bigl(\chi(-1)\bigr)^{2n-i(h^{-1})}$.
Note also that if $x=u(1-h)$ then
$$f_h(y,x)=f(y,-uh)=-f(-uh,y)=-f_{h^{-1}}(x,y),$$
showing that  $f_{h^{-1}}=-f_h^{\roman t}$, and hence that
$\delta_{h^{-1}}=\chi(-1)^{i(h)}\delta_h$. Thus we have shown that
$$\sum_{x\in I(g)\cap I(h^{-1})}\exp\gamma_{g,h}(x,x)
=\theta^{i(g)+i(h)-i(gh)}\delta_g\delta_h\delta_{gh},$$
and using Theorem 5.4 we deduce that
$$|I(g)|\,|I(h)|\,|I(gh)|^{-1}\sigma(g,h)
=\theta^{i(g)+i(h)-i(gh)}\delta_g\delta_h\delta_{gh}.$$
We can now conclude this section by proving the following theorem.
\proclaim{6.5 Theorem}If $f$ is a nondegenerate alternating form
$V\times V\to\F_p$, where $p$ is an odd prime, then the function $\mu$ defined 
on $G(f)$ by $\mu(g)=|I(g)|^{-1}\theta^{i(g)}\delta_g$ splits the factor set
$\sigma$. Furthermore, if $K'$ has an automorphism $\eta'$ inverting 
$\epsilon$ then $\mu(g)=\mu(g^{-1})^{\eta'}$ for all~$g\in G(f)$.\endproclaim
\demo{Proof} Define a factor set $\sigma'$ on $G=G(f)$ by
$$\sigma'(g,h)=\mu(g)^{-1}\mu(h)^{-1}\mu(gh)\sigma(g,h).$$
Our aim is to prove that $\sigma'(g,h)=1$ for all $g$ and $h$, and
our calculations above have established this whenever $(1-gh)$ is 
invertible. It is true whenever $h=1$, for then $\mu(h)=1$, and so
$\sigma'(g,1)=\sigma(g,1)=1$ . Let us check also
that it also holds whenever $h=g^{-1}$.

Since the form $\gamma_{g,g^{-1}}$ is zero and $I(g)=I(g^{-1})$, Theorem~5.4
gives $\sigma(g,g^{-1})=|I(g)|^{-1}$. It is clear that $\mu(1)=1$, and
$$\mu(g)\mu(g^{-1})=|I(g)|^{-2}\theta^{2i(g)}\delta_g\delta_{g^{-1}}=
|I(g)|^{-1}$$in view of the formulas for $\delta_{g^{-1}}$ and 
$\theta^2$. Hence $\sigma'(g,g^{-1})=1$.

Suppose now that $\psi$ is any projective representation with factor
set $\sigma'$. Then $\psi(1)$ is the identity (since $\sigma'(1,1)=1$), 
and $\psi(g)\psi(g^{-1})=\psi(1)$ since $\sigma'(g,g^{-1})=1$. 
Similarly, we have $\psi(g)\psi(h)=\psi(gh)$ whenever $1-gh$ is 
invertible. Writing $g_1=g^{-1}$ and $g_2=gh$ this gives
$\psi(g_1)\psi(g_2)=\psi(g_1g_2)$ whenever $1-g_2$ is invertible.
An obvious induction yields the same result whenever $g_2$ is a
product of elements $g$ with $1-g$ invertible. Since it is easily shown 
that the set of all such elements generates the whole symplectic 
group~$G$, we conclude that $\sigma'=1$.

Finally, observe that 
$$\mu(g^{-1})^{\eta'}=|I(g)|(\theta^{i(g)})^{\eta'}\delta_{g^{-1}}
=|I(g)|(\chi(-1)\theta)^{i(g)}\chi(-1)^{i(g)}\delta_g=\mu(g)$$
as required.
\enddemo

\heading 7. Towers of extraspecial groups\endheading

\noindent Suppose that $V$ is a $2n$-dimensional vector space over the field
$K=\F_q$ of characteristic~$p$. Let $\eta$ be an automorphism of $K$
satisfying $\eta^2=1$ and let $f$ be a $\eta$-sesquilinear form on~$V$.
More precisely, suppose that one of the following holds:
\roster
\item $\eta$ is the identity, $f$ is alternating and  $p\ne2$, or
\item $\eta$ is the identity, $V=W^*\oplus W$ and $f$ is the bilinear form
defined by $f((\alpha,x),(\beta,y))=x\beta$, or
\item $\eta$ is nontrivial and $f$ is $\eta$-Hermitian.
\endroster

Define the $\F_p$-bilinear form $\hat f=\hat f_\lambda$ by
$\hat f(x,y)=T(\lambda f(x,y))$, where $T\colon K\to\F_p$ is the trace map
and where $\lambda$ satisfies $\lambda^\eta\ne\lambda$ in case (3) and
$\lambda=1$ in cases (1) and (2). Then $\hat f-\hat f^{\roman t}$ is
nondegenerate and so $E(\hat f)$ is an extraspecial group of order
$p^{2n+1}$ whose isomorphism type may be determined by 4.1 and~4.2.

Let $\tilde\rho\colon E(\hat f)\to\gl(W)$ be a faithful absolutely
irreducible representation, where $W$ is a vector space over
a field $K'$ of characteristic $p'\ne p$. By the results of
Sections~5 and~6 we know that $\tilde\rho$ extends to a representation
$\bar\rho$ of $G(\hat f)\semiprod E(\hat f)$. If $E(\hat f)\cong D^n$ or
if $p$ is odd and $\roman{ord}_p(p')$ is odd we define $\rho$ to be the sum of
$\tilde\rho$ and $\tilde\rho^*$ (the contragredient of~$\tilde\rho$),
otherwise we define $\rho=\tilde\rho$. In the former case we define a
bilinear form $f'$ on $V'=W^*\oplus W$ as in Case~2 of Section~4, and note
that $\bar\rho^*\oplus\bar\rho$ is a representation of
$G(\hat f)\semiprod E(\hat f)$ extending $\rho$ and preserving~$f'$.
In the latter case we know by the results of Section~3 that $\rho$
preserves some $\eta'$-sesquilinear form~$f'$, and by the results of 
Sections~5 and~6 that the extension of $\rho$ also preserves~$f'$.
Thus in either case we have an embedding
$G(\hat f)\semiprod E(\hat f)\hookrightarrow G(f')$.


Let $f_1=\hat f$ and let $f_2=\widehat{f'}$ be defined in the same way,
but with $f'$ replacing~$f$. Since $G(f_2)=G(f')$ we have an embedding
$$G(f_1)\semiprod E(f_1)\hookrightarrow G(f_2),$$
giving rise to the group 
$$(G(f_1)\semiprod E(f_1))\semiprod E(f_2)$$
contained in $G(f_2)\semiprod E(f_2)$. Continuing this process, we may 
construct the iterated split extension
$$G(f_1)\semiprod E(f_1)\semiprod E(f_2)\semiprod\cdots\semiprod E(f_n).
\tag{$*$}$$

One group of this form may be constructed as follows. Let $f_i$  be a 
nondegenerate alternating form on a $2n_i$--dimensional vector space
over $\F_3$. By 4.1, $E(f_i)\cong E^{n_i}$ where $E$ is an extraspecial 
group of order 27 and exponent 3. There is a faithful irreducible
representation of $E(f_i)$ of degree $n_{i+1}=3^{n_i}$ over $\F_4$, which
preserves a nondegenerate Hermitian form $f_{i+1}$. By 4.2, $E(f_{i+1})
\cong Q^{n_{i+1}}$ and there is a 
$n_{i+2}=2^{n_{i+1}}$--dimensional representation of $E(f_{i+1})$ over 
$\F_3$ that preserves a nondegenerate bilinear form $f_{i+2}$, which is
alternating as $n_{i+1}$ is odd.
This is similar to the situation first considered except that
$f_i$ and $n_i$ are replaced by $f_{i+2}$ and $n_{i+2}/2$. If $f_1$
is an alternating form on a two--dimensional vector space over $\F_3$, then
$G(f_1)\cong {\text Sp}_2(\F_3)$ and we may therefore construct an iterated
split extension of the form
$$\text{Sp}_2(\F_3)\semiprod E\semiprod Q^3\semiprod E^4\semiprod Q^{81}
\semiprod E^{2^{80}}\semiprod\cdots\ .$$

In Section 1 we alluded to an extraspecial tower having $\gl_2(\F_3)$ as 
a quotient. Since $\text{Sp}_2(\F_3)\cong\text{SL}_2(\F_3)$, we seek an
extension of the above extraspecial tower by a cyclic group of order 2.
The construction of the larger extraspecial towers requires two steps.
If $G$ is a larger finite extraspecial tower, such as $S_3$ or $\gl_2(\F_3)$,
then we construct $G\semiprod E$ where $E$ is either an extraspecial 
2--group, or an extraspecial 3--group of the appropriate size.

First suppose that $G$ is one such larger extraspecial tower with a normal 
extraspecial 2--subgroup $E(f)\cong Q^n$ where $f$ is an $\eta$--Hermitian 
form defined on an $n$--dimensional vector space $V\cong E(f)\big/Z(E(f))$ 
over $\F_4$.
Suppose additionally that every $g\in G$ induces an $\alpha(g)$--semilinear
transformation $\bar g$ of $V$ such that
$$f(x{\bar g},y{\bar g})=f(x,y)^{\alpha(g)}\qquad\text{for all}\ x,y\in V,$$
where $\alpha(g)=\eta$ if $g\not\in G'$ and $\alpha(g)=1$ otherwise.
(This is the case for example if $G=\gl_2(\F_3)$.) If $g\not\in G'$, then 
$\bar g$ does not preserve $f$ (or even ${\hat f}_\lambda$), however, it
does preserve the quadratic form $V\rightarrow \F_2$ defined by
$x\mapsto {\hat f}_\lambda(x,x)=f(x,x)$. If 
$\rho: Q^n\rightarrow\gl_{2^n}(\F_3)$ is an irreducible representation that
preserves an alternating form $f'$, then by [Theorem 7, \glasby ] 
there is a representation $\tilde\rho$ extending $\rho$ that preserves 
$f'$ up to a sign. Hence we may construct the split extension
$G\semiprod E(f')$.

For the second inductive step, suppose that $G$ is an extraspecial tower with
a normal extraspecial 3--group $E(f)\cong E^n$ defined by an alternating form
$f$ on a $2n$--dimensional vector space $V$ over $\F_3$. Suppose additionally
that each $g\in G$ induces a linear transformation $\bar g$ of $V$ such that 
$$f(x{\bar g},y{\bar g})=\alpha(g)f(x,y)\qquad\text{for all}\ x,y\in V,$$
where $\alpha(g)=-1$ if $g\not\in G'$ and $\alpha(g)=1$ otherwise. 
An irreducible representation $\rho:E(f)\rightarrow\gl_{3^n}(\F_4)$
necessarily preserves some $\eta'$--Hermitian form $f'$, and by Section 6,
there is an extension $\tilde\rho$ of $\rho$ to $G'$ which also preserves $f'$.
Since every $g\not\in G'$ inverts $Z(E(f))$, the representation 
${\tilde\rho}^{(g)}:h\mapsto{\tilde\rho}(g^{-1}hg)$ of $G'$ is equivalent
to the representation 
${\tilde\rho}^{\eta'}:h\mapsto{\tilde\rho}(h)^{\eta'}$.
In this situation it is possible to extend $\tilde\rho$ to a crossed
representation of $G$. (Recall that $\sigma:G\rightarrow
\gl_m(\F)$ is a crossed representation if 
$\sigma(gh)=\sigma(g)^{\alpha(h)}\sigma(h)$ for
all $g,h\in G$ where $\alpha:G\rightarrow \text{Aut}(\F)$ is a 
homomorphism.) This crossed representation may be viewed as a representation
$G\rightarrow\gl_{2\cdot 3^n}(\F_2)$, and so we may construct the split
extension $G\semiprod E(f')$. This justifies the existence of the larger
extraspecial towers. These extraspecial towers have the property that
each of their normal subgroups are terms of their derived series. Furthermore,
some extraspecial towers provide examples of ``small'' soluble groups with
``large'' derived lengths (see [\gla ]) such as 
$$\gl_2(\F_3)\semiprod E\semiprod Q^3\semiprod E^4$$
which has order $2^{11}3^{13}$ and derived length 10.

The extraspecial towers discussed above were constructed using alternating and
Hermitian forms, one after the other in succession. There are many other 
possibilities of course. For example, each $f_i$ could be a Hermitian form 
acting on an $n_i$--dimensional vector space over $\F_{p_i^2}$ where $p_i$ 
is an odd prime. Provided $p_i|(p_{i+1}+1)$ for each $i$, we may construct 
the extraspecial tower $(*)$ in which $n_{i+1}=p_i^{n_i}$.


\Refs

\ref\no 1\by B. Bolt, T.G. Room and G.E. Wall\yr 1961\paper On the Clifford
collineation, transform and similarity groups, I, II\jour J. Austral. Math. 
Soc.\vol 2\pages 60--96\endref

\ref\no 2\by P. G\'erardin\yr 1977\paper Weil representations associated
to finite fields\jour J. Algebra\vol 46\pages 54--101\endref

\ref\no 3\by S.P. Glasby\yr 1989\paper The composition and derived lengths
of a soluble group\jour J. Algebra\vol 20\pages 406--413\endref

\ref\no 4\by S.P. Glasby\paper On the faithful representations of $2_\epsilon^{1+2n}\cdot 
O_{2n}^\epsilon(2)$ and $4\circ 2^{1+2n}\cdot Sp_{2n}(2)$ of
degree $2^n$\jour J. Austral. Math. Soc. (to appear)
\endref

\ref\no 5\by D. Gorenstein\yr 1980\book Finite Groups\publ Chelsea\publaddr 
New York\endref

\ref\no 6\by B. Huppert\yr 1967\book Endliche Gruppen, I\publ Springer--Verlag
\publaddr Berlin\endref

\ref\no 7\by G.E. Wall\yr 1963\paper On the conjugacy classes in the unitary, 
symplectic and orthogonal groups\jour J. Austral. Math. Soc.\vol 3\pages 1--62
\endref

\ref\no 8\by H. N. Ward\yr 1972\paper Representations of symplectic groups
\jour J. Algebra\vol 20\pages 182--195\endref

\ref\no 9\by H. N. Ward\yr 1974\paper Quadratic residue codes and 
symplectic groups \jour J. Algebra\vol 29\pages 150--171\endref
                                                               
\endRefs
\enddocument